\documentclass[12pt, reqno]{amsart}
\usepackage{amsmath, amsthm, amscd, amsfonts, amssymb, graphicx, hyperref,color}
\usepackage{enumitem}
\newtheorem{theorem}{Theorem}[section]
\newtheorem{lemma}[theorem]{Lemma}
\newtheorem{proposition}[theorem]{Proposition}

\theoremstyle{definition}

\theoremstyle{remark}
\newtheorem{remark}[theorem]{Remark}
\numberwithin{equation}{section}

\begin{document}
\setcounter{page}{1}


\title[On the circular numerical range of 5-by-5 partial isometries]
      {On the circular numerical range of 5-by-5 partial isometries}
       
\author[M. Naimi and M. Benharrat]{Mehdi Naimi$^{1}$, Mohammed Benharrat$^{2*}$}
\address{$^{1,2}$ Laboratoire de Math\'{e}matiques Fondamentales et Appliqu\'{e}es d'Oran. 
	D\'{e}partement de Formation Pr\'{e}paratoire en Sciences et Technology,
	Ecole Nationale Polytechnique d'Oran-Maurice Audin (Ex. ENSET d'Oran), 
	BP 1523 Oran-El M'naouar, 31000 Oran, Alg\'{e}rie.}
\email{\textcolor[rgb]{0.00,0.00,0.84}{mohammed.benharrat@enp-oran.dz, mohammed.benharrat@gmail.com}}
\email{\textcolor[rgb]{0.00,0.00,0.84}{naimimehdi.nm@gmail.com}}




\date{21/12/2020.
\newline \indent $^{*}$ Corresponding author\\
This work was supported by the Algerian research project: PRFU, no. C00L03ES310120180002.}
\begin{abstract}
  We prove, in some cases in term of kippenhahn curve,  that if 5-by-5 partial isometry whose numerical range  is a circular disc  then  its center is must be the origin. This gives a partial  affirmative answer of the  Conjecture 5.1. of [H. l. Gau et al., Linear and Multilinear Algebra, 64 (1) 2016, 14--35.], for the five dimensional case. 
\end{abstract}
 \maketitle

\section{introduction}
         
         
         Let $A$ be an $n\times n$ complex matrix, its numerical range $W(A)$ is, by definition the set of          
     complex numbers $$W(A)=\left\lbrace \langle Ax,x\rangle: \quad x\in\mathbb{C}^n, \quad \Vert x\Vert=1 \right\rbrace.$$
     
      It is well known  that $W(A)$ is a nonempty compact convex subset of $\mathbb{C}$, which contains all the eigenvalues of $A$ and therefore its convex hull. The matrix $A$ is said to be a partial isometry if
      it is isometric on the orthogonal complement of the kernel of $A$, $Ker(A)$. Assume that $A$ is a partial isometry  whose
     numerical range $ W(A)$ is a circular disc. The question  is whether the centre of $W(A)$ must be the origin.   Gau et al. ,\cite{GauPI}, gave an affirmative response  in case the dimension of the underlying
     space is at most $4$, as follows
    \begin{theorem}\cite[Theorem 2.1]{GauPI}\label{thm4by4}
    	If $A$ is an n-by-n ($n \leq 4$) partial isometry with $W(A)=\left\lbrace z\in \mathbb{C}:  \left|z-a \right|\leq r \right\rbrace$, ($r>0)$. then $a=0$.
    \end{theorem}
    
    Also have conjectured that  theorem remains valid if $A$ is  is an n-by-n partial isometry, see \cite[Conjecture 5.1.]{GauPI}. 
     
     We follow the same procedure as in \cite{GauPI},  we give an affirmative answer for this conjecture for 5-by-5 partial isometry in some cases in term of kippenhahn curve, to be more specific,  
         \begin{theorem}\label{5by5iso}
Let $A$ be  a $5\times 5$ partial isometry  matrix  $W(A)=\left\lbrace z\in \mathbb{C}:  \left|z-a \right|\leq r \right\rbrace$, ($r>0)$.
if the  kippenhahn curve $C_R(A)$ have one of the following  shape
\begin{enumerate}[label=(\roman*)]
\item $C_R(A)$ consist of three point and an elliptic disc.
\item $C_R(A)$ consist of two elliptic disc and a point.
\item $C_R(A)$ consist of  a curve of degree 4 with double tangent and an elliptic disc.
\end{enumerate}
 then $a=0$.
\end{theorem}
      
     Two successful approaches to establishing this result are  a canonical decomposition of $n \times n$ partial isometry matrix  and the  Kippenhahn's result for the numerical range of $n \times n$  matrix.
     
      It well  known that the numerical range of an $n \times n$  matrix $A$ is
     completely determined by its Kippenhahn polynomial
     $p_A(x, y, z) = det(x\Re( A )+ y\Im( A) + zI_n)$, where $\Re(A)=(A+A^*)/ 2 $ and $\Im(A)=(A-A^*)/ 2i$ are the real and the imaginary part of $A$, respectively,
     and $I_n$ denotes the $n \times n$  identity matrix. Let $C(A)$ be the dual of the algebraic curve defined to be the zero set of $p_A(x, y, z) = 0$, on the complex projective plane $\mathbb{CP}^2$, which consists of all equivalence classes of points in
     $\mathbb{C}^3\setminus(0.0.0)$ under the equivalence relation $\sim$, this relation is defined by
     $(x, y, z) \sim (x', y', z')$ 
     if and only if there is a nonzero $\lambda \in \mathbb{C}$ such that $(x, y, z) = \lambda(x', y', z')$.    Kippenhahn showed that  $W(A)$  is the convex hull of the real points of $C(A)$,
     see \cite{Kippenhahn1951} and its English translation \cite{Kippenhahn2008} for a
     detailed discussion of the connections between the polynomial $P_A$, and the
     numerical range of $A$. This characterization is used by  many authors  to answer the question when the numerical range of
     a matrix is an elliptic disc.  For $2\times 2$ matrices  a complete description of the numerical range  is
     well known, that is  $W(A)$ elliptic disk (with possibly degenerate interior), see \cite{Gustafson1997}. In \cite{Kippenhahn1951} Kippenhahn  showed
     that there are four classes of shapes which the numerical range of a  $3\times 3$ matrices 
     $A$ can assume. This was improved in \cite{Keeler1997} by expressing the
     conditions in terms of the eigenvalues and entries of $A$, which are easier to apply. By the same procedure,  these results are geralized for  $4\times 4$ matrices. Let us  mention here, that numerous results are known in this direction only for
     some special classes of matrices, as it happens, for partial isomertry, nilpotent, doubly stochastic  matrices (etc...), see \cite{stoc}, \cite{nil}, \cite{Benharrat}. But no unifying and general theory is yet available. 
     
     In this paper, firstly, with a similar approach to that used in \cite{Gau2006} and \cite{Datidar}, we will  give necessary and sufficient  condition for   $ 5\times 5$ matrix $A$ to have  an elliptical disc in her real Kippenhahn curve  $C_R(A)$. We also express those conditions in terms of eigenvalues and entries of $A$. All these conditions will be useful to construct a $ 5\times 5$ matrix with an elliptic numerical range (Section 2.). Secondly, we establish the result of Theorem \ref{5by5iso} in the special case of $S_5$-matrices  (  $A$ is $S_n$-matrix if $A$  is a contraction, the eigenvalues  of $A$ are all in the open unit disc $\mathbb{D}$ and
     $rank (I_n -A^*A) = 1$) (Section 3.). Finally, using the results of the two preceding section, we give the proof of   Theorem \ref{5by5iso} (Section 4.).
\section{Necessary and sufficient conditions for $ 5\times 5$ matrix  to have  an elliptical disc in its real Kippenhahn curve}
Let $A$ be a $ 5\times 5$ complex matrix.  Here we give necessary and sufficient conditions for which the associated curve $C_R(A)$ contains an elliptic disc or discs. It well known  that, by Schur's theorem, every square matrix is unitarily equivalent to an upper triangular matrix. So without loss of generality we can assume that 
\begin{equation}\label{upper-trian}
 A=\left[ \begin {array}{ccccc} \lambda_1&a_{12}&a_{13}&a_{14}&a_{15}
\\ \noalign{\medskip}0&\lambda_2&a_{23}&a_{24}&a_{25}
\\ \noalign{\medskip}0&0&\lambda_3&a_{34}&a_{35}
 \\ \noalign{\medskip}0&0&0&\lambda_4&a_{45}
 \\ \noalign{\medskip}0&0&0&0&\lambda_5
\end {array} \right].
\end{equation}  
and $\lambda_i=\alpha_i+i\beta_i$ , where $\alpha_j$ and $\beta_j$ are real for $j = 1, 2, 3, 4, 5$. 

Then, we have
\begin{align*}
P_A(x,y,z)&=det(x\Re( A )+ y\Im( A) + zI_5)\\
&=\left| \begin {array}{ccccc} \alpha_1x+\beta_1y+z&\frac{a_{12}}{2}(x-iy)&\frac{a_{13}}{2}(x-iy)&\frac{a_{14}}{2}(x-iy)&\frac{a_{15}}{2}(x-iy)
\\ \noalign{\medskip}(x+iy)\frac{\overline{a_{12}}}{2}&\alpha_2x+\beta_2y+z&\frac{a_{23}}{2}(x-iy)&\frac{a_{24}}{2}(x-iy)&\frac{a_{25}}{2}(x-iy)
\\ \noalign{\medskip}(x+iy)\frac{\overline{a_{13}}}{2}&(x+iy)\frac{\overline{a_{23}}}{2}&\alpha_3x+\beta_3y+z&\frac{a_{34}}{2}(x-iy)&\frac{a_{35}}{2}(x-iy)
 \\ \noalign{\medskip}(x+iy)\frac{\overline{a_{14}}}{2}&(x+iy)\frac{\overline{a_{24}}}{2}&(x+iy)\frac{\overline{a_{34}}}{2}&\alpha_4x+\beta_4y+z&\frac{a_{45}}{2}(x-iy)
 \\ \noalign{\medskip}(x+iy)\frac{\overline{a_{15}}}{2}&(x+iy)\frac{\overline{a_{25}}}{2}&(x+iy)\frac{\overline{a_{35}}}{2}&(x+iy)\frac{\overline{a_{45}}}{2}&\alpha_5x+\beta_5y+z
\end {array} \right|
\end{align*}
By straightforward calculus, we obtain
\begin{equation}\label{PAxyz}
P_A(x,y,z)= \prod\limits_{i=1}^5 (\alpha_ix+\beta_iy+z) -\dfrac{x^2+y^2}{4}Q(x,y,z),
\end{equation} 
where
\begin{align*}
Q(x,y,z)&=\sum\limits_{S_{ijk}S_{lm}}|a_{lm}|^2(\alpha_ix+\beta_iy+z)(\alpha_jx+\beta_jy+z)(\alpha_kx+\beta_ky+z)\\
&-\sum\limits_{S_{ij}S_{klm}}[x\Re(a_{kl}a_{lm}\overline{a_{km}})+y\Im (a_{kl}a_{lm}\overline{a_{km}})](\alpha_ix+\beta_iy+z)(\alpha_jx+\beta_jy+z)\\
&-\dfrac{x^2+y^2}{4}\sum\limits_{i=1}^{5}(\alpha_ix+\beta_iy+z)P_i
+ \sum\limits_{S_{i}S_{jklm}}(\alpha_ix+\beta_iy+z)[\dfrac{x^2-y^2}{2}\Re(a_{jk}a_{kl}a_{lm}\overline{a_{jm}})\\
&+xy\Im(a_{jk}a_{kl}a_{lm}\overline{a_{jm}})] 
+\dfrac{x^2+y^2}{4}\sum\limits_{S_{ij}S_{klm}}[x\Re(a_{kl}a_{lm}\overline{a_{km}})+y\Im (a_{kl}a_{lm}\overline{a_{km}})]|a_{ij}|^2 \\
&-\dfrac{x^2+y^2}{4}\sum\limits_{S_{ijkl}S_{iml}}[x\Re(a_{ij}a_{jk}a_{kl}\overline{a_{im}}\overline{a_{ml}})+y\Im (a_{ij}a_{jk}a_{kl}\overline{a_{im}}\overline{a_{ml}})]\\
&-\dfrac{x^2+y^2}{4}\sum\limits_{S_{ijk}S_{lm}S_{im}S_{lk}}[x\Re(a_{ij}a_{jk}a_{lm}\overline{a_{im}}\overline{a_{lk}})+y\Im(a_{ij}a_{jk}a_{lm}\overline{a_{im}}\overline{a_{lk}})]\\
&-\dfrac{1}{4}[(x^3-3xy^2)\Re(a_{12}a_{23}a_{34}a_{45}\overline{a_{15}}) +(-y^3+3yx^2)\Im(a_{12}a_{23}a_{34}a_{45}\overline{a_{15}})].
\end{align*}
With $S_{i_1i_2...i_n}$, for some $n\leq 5$,  denotes the collection of all n-tuples $( i_1,i_2,...,i_n)$ of natural numbers such that
$1\leq i_1< i_2< ...\leq i_5$ and
$$P_i=\sum\limits_{S_{jkm}S_{lm}}|a_{jk}|^2|a_{lm}|^2-\sum\limits_{S_{jkl}S_{jml}}\Re(a_{jk}a_{kl}\overline{a_{jm}}\overline{a_{ml}})-\sum\limits_{S_{jk}S_{lkm}}\Re(a_{jk}a_{lm}\overline{a_{jm}}\overline{a_{lk}}) $$
 for every $i=1,..,5$. The summation is taken with two by two different indexes witch means that $i\neq j\neq k \neq l\neq m$.

We begin by the following Lemma.
\begin{lemma}\label{2ellipse}
     Let A be a $5\times5$ matrix. Then the Kippenhahn curve $C_R(A)$
consists of two ellipses, one with foci $\lambda_1$, $\lambda_2$ and minor axis of length $r$, the other  with foci $\lambda_3$, $\lambda_4$ and minor axis of length $s$, and $\lambda_5$ if and only if 
 
  $P_A(x, y, z)=[(\alpha_1x+\beta_1y+z)(\alpha_2x+\beta_2y+z)-\frac{r^2}{4}(x^2+y^2)] [(\alpha_3x+\beta_3y+z)(\alpha_4x+\beta_4y+z)-\frac{s^2}{4}(x^2+y^2)](\alpha_5x+\beta_5y+z),$
  
  where $\lambda_j=\alpha_j+i\beta_j ,j=1,2,3,4,5$  and the $\alpha_j's$ and $\beta_j's$ are real.  
\end{lemma}

\begin{proof}
 Let $B=\begin{bmatrix}\lambda_1&r\\ 0&\lambda_2  \end{bmatrix} 
             \oplus 
 \begin{bmatrix}\lambda_3&s \\ 0&\lambda_4 \end{bmatrix}
             \oplus  \lambda_5$. As $C_R(A)=C_R(B)$, by duality the polynomials $P_A$ and $P_B$ have to be the same, therefore
             
  $P_A(x,y,z)=[(\alpha_1x+\beta_1y+z)(\alpha_2x+\beta_2y+z)-\frac{r^2(x^2+y^2)}{4}][(\alpha_3x+\beta_3y+z)(\alpha_4x+\beta_4y+z)-\frac{s^2(x^2+y^2)}{4}](\alpha_5x+\beta_5y+z).$
  
 The converse is clear.
\end{proof}
With the above lemma, we have the following theorem.

\begin{theorem}\label{th2ellipse}
Let's $A$ be in upper-triangular form \eqref{upper-trian}, Then the Kippenhahn curve $C_R(A)$
consists of two ellipses, one with foci $\lambda_p$, $\lambda_q$ and minor axis of length $r$, the other  with foci $\lambda_t$, $\lambda_v$ and minor axis of length $s$, and a point  $\lambda_w$ if and only if 
\begin{enumerate}[label=(\alph*)]

\item $r^2+s^2=\sum\limits_{S_{lm}}|a_{lm}|^2.$

\item $r^2(\lambda_w+\lambda_t+\lambda_v)+s^2(\lambda_w+\lambda_q+\lambda_p)   \\ =      \sum\limits_{S_{ijk}S_{lm}}|a_{lm}|^2(\lambda_i+\lambda_j+\lambda_k)-\sum\limits_{S_{klm}}a_{kl}a_{lm}\overline{a_{km}}.$ 

\item $r^2(\lambda_w\lambda_v + \lambda_w\lambda_t+\lambda_t\lambda_v)
+s^2 (\lambda_w\lambda_p + \lambda_w\lambda_q+\lambda_p\lambda_q)\\
 = \sum\limits_{S_{ijk}S_{lm}}|a_{lm}|^2(\lambda_i\lambda_j + \lambda_i\lambda_k+\lambda_j\lambda_k)-\sum\limits_{S_{ij}S_{klm}}(\lambda_i+\lambda_j)a_{kl}a_{lm}\overline{a_{km}}+\sum\limits_{S_{jklm}}a_{jk}a_{kl}a_{lm}\overline{a_{jm}}.$

\item $r^2\lambda_w\lambda_t\lambda_v
+s^2\lambda_w\lambda_p\lambda_q\\ =
 \sum\limits_{S_{ijk}S_{lm}}|a_{lm}|^2\lambda_i\lambda_j\lambda_k -\sum\limits_{S_{ij}S_{klm}}\lambda_i\lambda_j a_{kl}a_{lm}\overline{a_{km}}+\sum\limits_{S_{i}S_{jklm}}a_{jk}a_{kl}a_{lm}\overline{a_{jm}}\lambda_i \\-
 a_{12}a_{23}a_{34}a_{45}\overline{a_{15}}.$

\item $r^2\alpha_w\alpha_t\alpha_v  +s^2\alpha_w\alpha_p\alpha_q -\dfrac{r^2s^2}{4}\alpha_w  \\
=\sum\limits_{S_{ijk}S_{lm}}|a_{lm}|^2 \alpha_i\alpha_j\alpha_k 
-\sum\limits_{S_{ij}S_{klm}}\Re(a_{kl}a_{lm}\overline{a_{km}})\alpha_i\alpha_j 
-\dfrac{1}{4}\sum\limits_{i=1}^5P_i\alpha_i 
+\dfrac{1}{2}\sum\limits_{S_{i}S_{jklm}}\Re(a_{jk}a_{kl}a_{lm}\overline{a_{jm}})\alpha_i\\
+\dfrac{1}{4}\sum\limits_{S_{ij}S_{klm}}\Re(a_{kl}a_{lm}\overline{a_{km}}) |a_{ij}|^2 
-\dfrac{1}{4}\sum\limits_{S_{ijkl}S_{iml}}\Re(a_{ij}a_{jk}a_{kl}\overline{a_{im}}\overline{a_{ml}})\\
-\dfrac{1}{4}\sum\limits_{S_{ijk}S_{lm}S_{im}S_{lk}}\Re(a_{ij}a_{jk}a_{lm}\overline{a_{im}}\overline{a_{lk}})
-\dfrac{1}{4}\Re(a_{12}a_{23}a_{34}a_{45}\overline{a_{15}}).$

\item $r^2\beta_w\beta_t\beta_v  +s^2\beta_w\beta_p\beta_q -\dfrac{r^2s^2}{4}\beta_w \\
 =\sum\limits_{S_{ijk}S_{lm}}|a_{lm}|^2 \beta_i\beta_j\beta_k 
-\sum\limits_{S_{ij}S_{klm}}\Im(a_{kl}a_{lm}\overline{a_{km}})\beta_i\beta_j 
-\dfrac{1}{4}\sum\limits_{i=1}^5P_i\beta_i 
-\dfrac{1}{2}\sum\limits_{S_{i}S_{jklm}}\Re(a_{jk}a_{kl}a_{lm}\overline{a_{jm}})\beta_i\\
+\dfrac{1}{4}\sum\limits_{S_{ij}S_{klm}}\Im(a_{kl}a_{lm}\overline{a_{km}}) |a_{ij}|^2 
-\dfrac{1}{4}\sum\limits_{S_{ijkl}S_{iml}}\Im(a_{ij}a_{jk}a_{kl}\overline{a_{im}}\overline{a_{ml}})\\
-\dfrac{1}{4}\sum\limits_{S_{ijk}S_{lm}S_{im}S_{lk}}\Im(a_{ij}a_{jk}a_{lm}\overline{a_{im}}\overline{a_{lk}})
+\dfrac{1}{4}\Im(a_{12}a_{23}a_{34}a_{45}\overline{a_{15}}).$

 \item  $r^2  (\alpha_w\alpha_t   + \alpha_w\alpha_v  +\alpha_v\alpha_t ) 
 +s^2 (\alpha_w\alpha_p   + \alpha_w\alpha_q  +\alpha_p\alpha_q ) -\dfrac{r^2s^2}{4} \\
 =   \sum\limits_{S_{ijk}S_{lm}}|a_{lm}|^2 (\alpha_i\alpha_j+\alpha_i\alpha_k + \alpha_j\alpha_k)
- \sum\limits_{S_{ij}S_{klm}}\Re(a_{kl}a_{lm}\overline{a_{km}})(\alpha_i+\alpha_j)\\
-\dfrac{1}{4}\sum\limits_{i=1}^5P_i
+ \dfrac{1}{2}\sum\limits_{S_{jklm}}\Re(a_{jk}a_{kl}a_{lm}\overline{a_{jm}}).
 $ 
\end{enumerate}
where 
$P_i=\sum\limits_{S_{jkm}S_{lm}}|a_{jk}|^2|a_{lm}|^2-\sum\limits_{S_{jkl}S_{jml}}\Re(a_{jk}a_{kl}\overline{a_{jm}}\overline{a_{ml}})-\sum\limits_{S_{jk}S_{lkm}}\Re(a_{jk}a_{lm}\overline{a_{jm}}\overline{a_{lk}}).$\\ for every $i=1,..,5$. 
\end{theorem}
\begin{proof}
By lemma \ref{2ellipse} we get 
\begin{align*}
 P_A(x, y, z) &= (\alpha_wx+\beta_wy+z)(\alpha_px+\beta_py+z)(\alpha_qx+\beta_qy+z)(\alpha_tx+\beta_ty+z)(\alpha_vx+\beta_vy+z)\\
 & -\dfrac{x^2+y^2}{4}[r^2(\alpha_wx+\beta_wy+z)(\alpha_tx+\beta_ty+z)(\alpha_vx+\beta_vy+z)\\ &+s^2(\alpha_wx+\beta_wy+z)(\alpha_px+\beta_py+z)(\alpha_qx+\beta_qy+z)
 \\ &-\dfrac{x^2+y^2}{4}r^2s^2(\alpha_wx+\beta_wy+z)].
\end{align*}

Comparing the previous formula of $P_A(x, y, z)$ with \eqref{PAxyz}. We find 
   \begin{align*}
 Q(x, y, z) &=
 r^2(\alpha_wx+\beta_wy+z)(\alpha_tx+\beta_ty+z)(\alpha_vx+\beta_vy+z)\\ &+s^2(\alpha_wx+\beta_wy+z)(\alpha_px+\beta_py+z)(\alpha_qx+\beta_qy+z)
 \\ &-\dfrac{x^2+y^2}{4}r^2s^2(\alpha_wx+\beta_wy+z).
 \end{align*}
Computing the coefficients of $x^3,y^3,z^3,x^2y,xy^2,x^2z,xz^2,y^2z,yz^2,xyz$ by identification, we find, respectively 
\begin{enumerate}
\item $r^2\alpha_w\alpha_t\alpha_v  +s^2\alpha_w\alpha_p\alpha_q -\dfrac{r^2s^2}{4}\alpha_w  \\
=\sum\limits_{S_{ijk}S_{lm}}|a_{lm}|^2 \alpha_i\alpha_j\alpha_k 
-\sum\limits_{S_{ij}S_{klm}}\Re(a_{kl}a_{lm}\overline{a_{km}})\alpha_i\alpha_j 
-\dfrac{1}{4}\sum\limits_{i=1}^5P_i\alpha_i 
+\dfrac{1}{2}\sum\limits_{S_{i}S_{jklm}}\Re(a_{jk}a_{kl}a_{lm}\overline{a_{jm}})\alpha_i\\
+\dfrac{1}{4}\sum\limits_{S_{ij}S_{klm}}\Re(a_{kl}a_{lm}\overline{a_{km}}) |a_{ij}|^2 
-\dfrac{1}{4}\sum\limits_{S_{ijkl}S_{iml}}\Re(a_{ij}a_{jk}a_{kl}\overline{a_{im}}\overline{a_{ml}})\\
-\dfrac{1}{4}\sum\limits_{S_{ijk}S_{lm}S_{im}S_{lk}}\Re(a_{ij}a_{jk}a_{lm}\overline{a_{im}}\overline{a_{lk}})
-\dfrac{1}{4}\Re(a_{12}a_{23}a_{34}a_{45}\overline{a_{15}}).$

\item $r^2\beta_w\beta_t\beta_v  +s^2\beta_w\beta_p\beta_q -\dfrac{r^2s^2}{4}\beta_w \\
 =\sum\limits_{S_{ijk}S_{lm}}|a_{lm}|^2 \beta_i\beta_j\beta_k 
-\sum\limits_{S_{ij}S_{klm}}\Im(a_{kl}a_{lm}\overline{a_{km}})\beta_i\beta_j 
-\dfrac{1}{4}\sum\limits_{i=1}^5P_i\beta_i 
-\dfrac{1}{2}\sum\limits_{S_{i}S_{jklm}}\Re(a_{jk}a_{kl}a_{lm}\overline{a_{jm}})\beta_i\\
+\dfrac{1}{4}\sum\limits_{S_{ij}S_{klm}}\Im(a_{kl}a_{lm}\overline{a_{km}}) |a_{ij}|^2 
-\dfrac{1}{4}\sum\limits_{S_{ijkl}S_{iml}}\Im(a_{ij}a_{jk}a_{kl}\overline{a_{im}}\overline{a_{ml}})\\
-\dfrac{1}{4}\sum\limits_{S_{ijk}S_{lm}S_{im}S_{lk}}\Im(a_{ij}a_{jk}a_{lm}\overline{a_{im}}\overline{a_{lk}})
+\dfrac{1}{4}\Im(a_{12}a_{23}a_{34}a_{45}\overline{a_{15}}).$

 \item  $r^2 +s^2 = \sum\limits_{S{lm}}|a_{lm}|^2.$

 \item  $r^2 (\beta_w\alpha_t\alpha_v   +\beta_t\alpha_w\alpha_v  +\beta_v\alpha_w\alpha_t ) 
 +  s^2 (\beta_w\alpha_p\alpha_q  +\beta_p\alpha_w\alpha_q  +\beta_q\alpha_w\alpha_p )-\dfrac{r^2s^2}{4}\beta_w\\
 =\sum\limits_{S_{ijk}S_{lm}}|a_{lm}|^2 (\alpha_i\alpha_j\beta_k+\alpha_i\alpha_k\beta_j + \alpha_j\alpha_k\beta_i)
 -\sum\limits_{S_{ij}S_{klm}}\Im(a_{kl}a_{lm}\overline{a_{km}})\alpha_i\alpha_j \\
 -\sum\limits_{S_{ij}S_{klm}}\Re(a_{kl}a_{lm}\overline{a_{km}})(\alpha_i\beta_j +\alpha_j\beta_i) 
 -\dfrac{1}{4}\sum\limits_{i=1}^5P_i\beta_i 
 +\dfrac{1}{2}\sum\limits_{S_{i}S_{jklm}}\Re(a_{jk}a_{kl}a_{lm}\overline{a_{jm}})\beta_i\\
 +\sum\limits_{S_{i}S_{jklm}}\Im(a_{jk}a_{kl}a_{lm}\overline{a_{jm}})\alpha_i
 + \dfrac{1}{4}\sum\limits_{S_{ij}S_{klm}}\Im(a_{kl}a_{lm}\overline{a_{km}}) |a_{ij}|^2\\
 -\dfrac{1}{4}\sum\limits_{S_{ijkl}S_{iml}}\Im(a_{ij}a_{jk}a_{kl}\overline{a_{im}}\overline{a_{ml}})
 -\dfrac{1}{4}\sum\limits_{S_{ijk}S_{lm}S_{im}S_{lk}}\Im(a_{ij}a_{jk}a_{lm}\overline{a_{im}}\overline{a_{lk}}) \\
-\dfrac{3}{4}\Im(a_{12}a_{23}a_{34}a_{45}\overline{a_{15}})  
.   $

  \item   $r^2 (\alpha_w\beta_t\beta_v   +\alpha_t\beta_w\beta_v  +\alpha_v\beta_w\beta_t )
 +  s^2 (\alpha_w\beta_p\beta_q  +\alpha_p\beta_w\beta_q  +\alpha_q\beta_w\beta_p )-\dfrac{r^2s^2}{4}\alpha_w\\
 = \sum\limits_{S_{ijk}S_{lm}}|a_{lm}|^2 (\alpha_i\beta_j\beta_k+\alpha_j\beta_i\beta_k + \alpha_k\beta_i\beta_j)
 -\sum\limits_{S_{ij}S_{klm}}\Re(a_{kl}a_{lm}\overline{a_{km}})\beta_i\beta_j \\
 -\sum\limits_{S_{ij}S_{klm}}\Im(a_{kl}a_{lm}\overline{a_{km}})(\alpha_i\beta_j +\alpha_j\beta_i)  
 -\dfrac{1}{4}\sum\limits_{i=1}^5P_i\alpha_i  
 -\dfrac{1}{2}\sum\limits_{S_{i}S_{jklm}}\Re(a_{jk}a_{kl}a_{lm}\overline{a_{jm}})\alpha_i\\ 
 +\sum\limits_{S_{i}S_{jklm}}\Im(a_{jk}a_{kl}a_{lm}\overline{a_{jm}})\beta_i 
 + \dfrac{1}{4}\sum\limits_{S_{ij}S_{klm}}\Re(a_{kl}a_{lm}\overline{a_{km}}) |a_{ij}|^2\\ 
 -\dfrac{1}{4}\sum\limits_{S_{ijkl}S_{iml}}\Re(a_{ij}a_{jk}a_{kl}\overline{a_{im}}\overline{a_{ml}}) 
 -\dfrac{1}{4}\sum\limits_{S_{ijk}S_{lm}S_{im}S_{lk}}\Re(a_{ij}a_{jk}a_{lm}\overline{a_{im}}\overline{a_{lk}}) \\
+\dfrac{3}{4}\Re(a_{12}a_{23}a_{34}a_{45}\overline{a_{15}}) .  $

 \item  $r^2  (\alpha_w\alpha_t   + \alpha_w\alpha_v  +\alpha_v\alpha_t ) 
 +s^2 (\alpha_w\alpha_p   + \alpha_w\alpha_q  +\alpha_p\alpha_q ) -\dfrac{r^2s^2}{4} \\
 =   \sum\limits_{S_{ijk}S_{lm}}|a_{lm}|^2 (\alpha_i\alpha_j+\alpha_i\alpha_k + \alpha_j\alpha_k)
- \sum\limits_{S_{ij}S_{klm}}\Re(a_{kl}a_{lm}\overline{a_{km}})(\alpha_i+\alpha_j)\\
-\dfrac{1}{4}\sum\limits_{i=1}^5P_i
+ \dfrac{1}{2}\sum\limits_{S_{jklm}}\Re(a_{jk}a_{kl}a_{lm}\overline{a_{jm}}).
 $

 \item $r^2(\alpha_w+\alpha_t+\alpha_v)+s^2(\alpha_w+\alpha_p+\alpha_q)\\
 =\sum\limits_{S_{ijk}S_{lm}}|a_{lm}|^2 (\alpha_i + \alpha_j + \alpha_k) 
 -\sum\limits_{S_{klm}}\Re(a_{kl}a_{lm}\overline{a_{km}}).
 $

  \item  $r^2 (\beta_w\beta_t   + \beta_w\beta_v  +\beta_v\beta_t ) 
 +s^2 (\beta_w\beta_p   + \beta_w\beta_q  +\beta_p\beta_q ) -\dfrac{r^2s^2}{4}\\
 =   \sum\limits_{S_{ijk}S_{lm}}|a_{lm}|^2 (\beta_i\beta_j+\beta_i\beta_k + \beta_j\beta_k)
- \sum\limits_{S_{ij}S_{klm}}\Im(a_{kl}a_{lm}\overline{a_{km}})(\beta_i+\beta_j)\\
-\dfrac{1}{4}\sum_{i=1}^5P_i
-\dfrac{1}{2}\sum_{S_{jklm}}\Re(a_{jk}a_{kl}a_{lm}\overline{a_{jm}}).$ 
 
 \item $r^2(\beta_w+\beta_t+\beta_v)+s^2(\beta_w+\beta_p+\beta_q)
 = \sum\limits_{S_{ijk}S_{lm}}\left| a_{lm}\right| ^2 (\beta_i + \beta_j + \beta_k) 
 -\sum_{S_{klm}}\Im(a_{kl}a_{lm}\overline{a_{km}}).
 $
 \item 
 $r^2(\alpha_w\beta_t +\alpha_w\beta_v+\alpha_t\beta_w+\alpha_t\beta_v+\alpha_v\beta_w+\alpha_v\beta_t)+ s^2(\alpha_w\beta_p +\alpha_w\beta_q+\alpha_p\beta_w+\alpha_p\beta_q+\alpha_q\beta_w+\alpha_q\beta_p)\\
  =\sum\limits_{S_{ijk}S_{lm}}|a_{lm}|^2 (\alpha_i\beta_j +\alpha_j\beta_i+\alpha_i\beta_k+\alpha_k\beta_i+\alpha_j\beta_k+\alpha_k\beta_j)- \sum\limits_{S_{ij}S_{klm}}\Re(a_{kl}a_{lm}\overline{a_{km}})(\beta_i+\beta_j) \\
  - \sum\limits_{S_{ij}S_{klm}}\Im(a_{kl}a_{lm}\overline{a_{km}})(\alpha_i+\alpha_j) 
+\sum\limits_{S_{jklm}}\Im(a_{jk}a_{kl}a_{lm}\overline{a_{jm}})$.\\

\end{enumerate}

   Note that the combination of $(1)$, $(2)$, $(4)$ and $(5)$ is equivalent to $(d)$, $(e)$ and $(f)$, since 
$(1)-(5)-i(2)+i(4)$ yields (d).

 $(6)$, $(8)$ and $(10)$ is equivalent to $(c)$ and $(g)$ as $(6)-(8) + i(10)$ yields $(c)$.
 
$(7)$ and $(9)$ is equivalent to $(b)$ as it's follows from $(7)+i(9)$.

 This completes the proof.
\end{proof} 
%
\begin{remark} 

\begin{itemize}
\item The last theorem was obtained in \cite[Theorem 2.2]{Datidar}, but, there exists a gap in this Theorem, we shall point out that the  equations $(4')$ and $(5')$ of the proof of \cite[Theorem 2.2, p 726 ]{Datidar} are missed. We give here corrected formulas.  Also, the authors claimed in the end of the proof that the combination of $(1')$, $(2')$ $(4')$ and $(5')$ is equivalent to $(d)$ (see page 727). But,  we should be  add the condition $(e)$, $(f)$, and $(g)$ to ensure the converse implication.
\item If we take $s=0$ in the previous theorem then $C_R(A)$ is an  ellipse and three points.
\end{itemize}
\end{remark}
%
%
In the following, another case comes materialized, it's this when $C_R(A)$ contain an ellipse and  a curve of degree 4 with a double tangent.  Using the same conditions derived in   \cite{Rodman}, to determine whether the numerical range  boundary  of a given $3\times 3$ irreducible  matrix $A$ has a flatness, these conditions are given in term of geometrical property of flatness.
  
Let $L$ be the supporting line of the convex set $W(A)$ containing the flatness and   perpendicular to the line which pass through the origin and forms angle $\theta$ from the positive $x-$axis, and let $\mu$ be the (signed) distance from
the origin to $L$. It is seen that $\mu$ is the largest eigenvalue of $\Re(e^{-i\theta}A)$  (cf.\cite{Rodman},\cite{Kippenhahn1951}).
%
%
\begin{lemma}\label{flat}
     Let A be a $5\times 5$ matrix. Then the Kippenhahn curve $C_R(A)$
consists of one ellipse with foci $\lambda_1$,$\lambda_2$ and minor axis of length $r$, and a curve of degree 4 with a double tangent and foci at   $\lambda_3$ , $\lambda_4$ and $\lambda_5$ if and only if the following condition hold
\begin{enumerate}[label=(\roman*)]
\item  there exist $\theta \in [0,2\pi[$ and a  real $\mu$ such that 
    \begin{align*} 
  P_A(x, y, z)& =\left[(\alpha_1x+\beta_1y+z)(\alpha_2x+\beta_2y+z)-\frac{r^2}{4}(x^2 +y^2)
  \right]\\ & [(\alpha_3x+\beta_3y+z)(\alpha_4x+\beta_4y+z)(\alpha_5x+\beta_5y+z)\\ & -(\alpha_3x+\beta_3y+z)(x^2 +y^2)(\Re(e^{-i\theta}\lambda_4) +\mu)(\Re(e^{-i\theta}\lambda_5) +\mu)  \\ & 
  -(\alpha_4x+\beta_4y+z)(x^2 +y^2)(\Re(e^{-i\theta}\lambda_3) +\mu)(\Re(e^{-i\theta}\lambda_5) +\mu)  \\ & 
  -(\alpha_5x+\beta_5y+z)(x^2 +y^2)(\Re(e^{-i\theta}\lambda_3) +\mu)(\Re(e^{-i\theta}\lambda_4) +\mu)  \\ &+(x-iy)(x^2 +y^2)(\Re(e^{-i\theta}\lambda_3) +\mu)(\Re(e^{-i\theta}\lambda_4) +\mu)(\Re(e^{-i\theta}\lambda_5) +\mu)e^{i\theta} \\ & + (x+iy)(x^2 +y^2)(\Re(e^{-i\theta}\lambda_3) +\mu)(\Re(e^{-i\theta}\lambda_4) +\mu)(\Re(e^{-i\theta}\lambda_5) +\mu)e^{-i\theta}.
    \end{align*} 
\item  $(\Re(e^{-i\theta}\lambda_3) +\mu)(\Re(e^{-i\theta}\lambda_4) +\mu)(\Re(e^{-i\theta}\lambda_5) +\mu)\neq 0$.
\item $\lambda_j(\Re(e^{-i\theta}\lambda_k) +\mu)+ \lambda_k(\Re(e^{-i\theta}\lambda_j) +\mu) - 2(\Re(e^{-i\theta}\lambda_j) +\mu)(\Re(e^{-i\theta}\lambda_k) +\mu) e^{i\theta}\neq \lambda_i((\Re(e^{-i\theta}\lambda_k) +\mu) + (\Re(e^{-i\theta}\lambda_j) +\mu))$, for every  $3\leq i\neq j\neq k \leq 5.$
 \end{enumerate}
\end{lemma}
%
%
\begin{proof} 
     Let $B=\begin{bmatrix}\lambda_1&r
\\ 0&\lambda_2
   \end{bmatrix} \oplus C$ where $C$ is 3-by-3 irreducible matrix whose Kippenhahn curve is of degree 4
   with a double tangent, and eigenvalues at   $\lambda_3$ , $\lambda_4$ and $\lambda_5$.   \\
   Since $C_R(A) = C_R(B)$, the polynomials $P_A$ and $P_B$ have to be the same.\\
 $$p_B(x,y,z)=\left[(\alpha_1x+\beta_1y+z)(\alpha_2x+\beta_2y+z)-\frac{r^2}{4}(x^2 +y^2)
  \right]p_C(x,y,z)$$
 Put $C$ in an upper triangular form
  	$\begin{bmatrix}
	\lambda_3&a&b
	\\ 0&\lambda_4&c
	\\0&0&\lambda_5
	\end{bmatrix}$.\\
 As we assume that $\partial W(C)$ has a flat of  portion (containing in the supporting line $L$), let's $\theta \in [0,2\pi[$ be the angular between x-axis and the line witch pass through the origin and perpendicular to $L$, so $e^{-i\theta}C$ has a vertical flatness. According to Kippenhahn's 
 classification, $\Re( e^{-i\theta}C)$ must have a multiple eigenvalue, so there exist a real $\mu$ such that  	$$\Re( e^{-i\theta}C)+I \mu=\begin{bmatrix}
	\\ \Re( e^{-i\theta}\lambda_3)+\mu &e^{-i\theta}a/2&e^{-i\theta}b/2
	\\ e^{i\theta}\overline{a}/2   &    \Re( e^{-i\theta}\lambda_4)+\mu  &     e^{-i\theta}c/2
	\\e^{i\theta}\overline{b}/2   &     e^{i\theta}\overline{c}/2   &    \Re( e^{-i\theta}\lambda_5)+\mu
	\end{bmatrix}$$
 has rank one, because if otherwise  $\Re( e^{-i\theta}C)+I\mu$ has zero rank, then $\Re( e^{-i\theta}C)$ and $\Im( e^{-i\theta}C)$ commute, and $C$ is therefore reducible, while due to the latter all $2\times 2$ minors of $\Re( e^{-i\theta}C)+I\mu$ are equal to zero, witch gives
\begin{align}\label{eqRodman}\begin{split}
| a |^2 &= 4(\Re(e^{-i\theta}\lambda_3) +\mu)(\Re(e^{-i\theta}\lambda_4) +\mu).\\
| b |^2 &= 4(\Re(e^{-i\theta}\lambda_3) +\mu)(\Re(e^{-i\theta}\lambda_5) +\mu).\\
| c |^2 &= 4(\Re(e^{-i\theta}\lambda_4) +\mu)(\Re(e^{-i\theta}\lambda_5) +\mu).\\
a\overline{b}&=2(\Re(e^{-i\theta}\lambda_3) + \mu)\overline{c}e^{i\theta}.\end{split}
\end{align}
              It's easy to see from  equations above that if one of off-diagonal $a,b$ or $c$ 
              is zero then at least two of them are equal to zero, contradiction
              with irreducibility of $C$, so $abc\neq 0$. \\
               On the other hand              
\begin{align*}
p_C(x,y,z)&=det\begin{bmatrix}
	(\alpha_3x+\beta_3y+z)&a/2(x-iy)&b/2(x-iy)
	\\ \overline{a}/2(x+iy) &(\alpha_4x+\beta_4y+z)&c/2(x-iy)
	\\\overline{b}/2(x+iy)&\overline{c}/2(x+iy)&(\alpha_5x+\beta_5y+z)
	\end{bmatrix} \\
	&=(\alpha_3x+\beta_3y+z)(\alpha_4x+\beta_4y+z)(\alpha_5x+\beta_5y+z)\\
&-\dfrac{1}{4}(x^2+y^2)[(\alpha_3x+\beta_3y+z)|c|^2+(\alpha_4x+\beta_4y+z)|b|^2+(\alpha_5x+\beta_5y+z)|a|^2]\\
&+\dfrac{1}{8}(x^2+y^2)(x-iy)a\overline{b}c+\dfrac{1}{8}(x^2+y^2)(x+iy)\overline{a}\overline{c}b.
\end{align*}
           Combining this relation with \eqref{eqRodman} and $a\overline{b}c\neq 0$ we get successively
           $(i)$ and $(ii)$,
 moreover the polynomial $P_C$ (and therefore the matrix $C$) is irreducible, consequently $P_C$ can not be factored into threes linear factors, or into a quadratic factors and a linear factor, we note that  linear factors in the left-hand side of the equation $P_C(x,y,z) = 0$ are corresponding always to eigenvalues of the matrix $C$ (c.f  \cite{Kippenhahn2008},\cite{Kippenhahn1951} and  \cite{Rodman}).\\
  Let's assume that $P_C(x,y,z)$ have linear factor $(\alpha_i x + \beta_i y+z)$ witch correspond to $\lambda_i:i=3,4,5.$
  and with 
\begin{align*}   
   P_C(x,y,z)&= (\alpha_3x+\beta_3y+z)(\alpha_4x+\beta_4y+z)(\alpha_5x+\beta_5y+z)\\ & -(\alpha_3x+\beta_3y+z)(x^2 +y^2)\mu_4\mu_5  \\ & 
  -(\alpha_4x+\beta_4y+z)(x^2 +y^2)\mu_3\mu_5  \\ & 
  -(\alpha_5x+\beta_5y+z)(x^2 +y^2)\mu_3\mu_4  \\
   &+(x-iy)(x^2 +y^2)\mu_3\mu_4\mu_5 e^{i\theta} \\
   & + (x+iy)(x^2 +y^2)\mu_3\mu_4\mu_5e^{-i\theta}.
    \end{align*}  
  where $\mu_i=\Re(e^{-i\theta}\lambda_i) +\mu$, 
 we can see that for two by two in-equal $i,j,k\in \lbrace 1,2,3\rbrace$, $(\alpha_i x + \beta_i y+z)$ devises 
 \begin{align*}   
&(\alpha_jx+\beta_jy+z)(x^2 +y^2)\mu_i\mu_k  \\
 &  +(\alpha_kx+\beta_ky+z)(x^2 +y^2)\mu_i\mu_j  \\ &
 -(x-iy)(x^2 +y^2)\mu_i\mu_j\mu_ke^{i\theta} \\
  & - (x+iy)(x^2 +y^2)\mu_i\mu_j\mu_ke^{-i\theta}.
    \end{align*}  
 in other words, the coefficients of $z,\alpha_ix,\beta_iy$ in the last polynomial are equals, witch gives 
 \begin{equation}\label{ir1}
 \beta_i[\alpha_j\mu_i\mu_k + \alpha_k\mu_i\mu_j -2 \mu_i\mu_j\mu_k \cos(\theta)]= \alpha_i[\beta_j\mu_i\mu_k + \beta_k\mu_i\mu_j - 2\mu_i\mu_j\mu_k \sin(\theta)]
 \end{equation}
  \begin{equation}\label{ir2}
\alpha_j\mu_i\mu_k + \alpha_k\mu_i\mu_j -2 \mu_i\mu_j\mu_k \cos(\theta)= \alpha_i[\mu_i\mu_k + \mu_i\mu_j ]
 \end{equation}
 and
  \begin{equation}\label{ir3}
\beta_j\mu_i\mu_k + \beta_k\mu_i\mu_j -2 \mu_i\mu_j\mu_k \sin(\theta)= \beta_i[\mu_i\mu_k + \mu_i\mu_j]
 \end{equation}
One can see that \eqref{ir1},\eqref{ir2} and \eqref{ir3} are equivalents to 
$\lambda_j\mu_k + \lambda_k\mu_j - 2\mu_j\mu_k e^{i\theta}= \lambda_i(\mu_k + \mu_j)$ and  hence $(iii)$
 
the converse is clear.
\end{proof}
In what follows, by $\mu_j$ we note $\Re(e^{-i\theta}\lambda_j)+\mu$ for every $j=1,...,5$.
\begin{theorem}\label{th4degree}
     Let $A$ be in upper-triangular form \eqref{upper-trian}. Then $C_R(A)$ 
consists of one ellipse with foci $\lambda_q$,$\lambda_p$ and minor axis of length $r$, and a curve of degree 4 with a double tangent and foci at   $\lambda_w$ , $\lambda_t$ and $\lambda_v$ if and only if 
there exist $\theta \in [0,2\pi[$ and a  real $\mu$ such that 
\begin{enumerate}[label=(\alph*)]

\item $r^2+4(\mu_t\mu_v+\mu_w\mu_v+\mu_w\mu_t)=\sum\limits_{S_{lm}}|a_{lm}|^2.$

\item $r^2(\lambda_w+\lambda_t+\lambda_v)+ 4[
(\lambda_w\mu_t\mu_v+\lambda_t\mu_w\mu_v+\lambda_v\mu_t\mu_w)- 2\mu_w\mu_t\mu_ve^{i\theta} \\ 
+(\lambda_p+\lambda_q)(\mu_t\mu_v+\mu_w\mu_v+\mu_w\mu_t)]  \\   =   \sum\limits_{S_{ijk}S_{lm}}|a_{lm}|^2(\lambda_i+\lambda_j+\lambda_k)-\sum\limits_{S_{klm}}a_{kl}a_{lm}\overline{a_{km}}.$ 

\item $r^2(\lambda_w\lambda_t + \lambda_t\lambda_v+\lambda_w\lambda_v)+4[ \lambda_p\lambda_q(\mu_t\mu_v+\mu_w\mu_v+\mu_t\mu_w)\\ +(\lambda_p+\lambda_q)(\lambda_w\mu_t\mu_v+\lambda_t\mu_w\mu_v+\lambda_v\mu_t\mu_w)
-(\lambda_p+\lambda_q)(2\mu_w\mu_t\mu_ve^{i\theta})]\\ 
=\sum\limits_{S_{ijk}S_{lm}}|a_{lm}|^2(\lambda_i\lambda_j + \lambda_i\lambda_k+\lambda_j\lambda_k)-\sum\limits_{S_{ij}S_{klm}}(\lambda_i+\lambda_j)a_{kl}a_{lm}\overline{a_{km}}+\sum\limits_{S_{jklm}}a_{jk}a_{kl}a_{lm}\overline{a_{jm}}.$

\item $r^2(\lambda_w\lambda_t\lambda_v) +4[ \lambda_p\lambda_q((\lambda_w\mu_t\mu_v+\lambda_t\mu_w\mu_v+\lambda_v\mu_t\mu_w)- 2\mu_w\mu_t\mu_ve^{i\theta})]\\ =
 \sum\limits_{S_{ijk}S_{lm}}|a_{lm}|^2\lambda_i\lambda_j\lambda_k -\sum\limits_{S_{ij}S_{klm}}\lambda_i\lambda_j a_{kl}a_{lm}\overline{a_{km}}+\sum\limits_{S_{i}S_{jklm}}a_{jk}a_{kl}a_{lm}\overline{a_{jm}}\lambda_i \\-
 a_{12}a_{23}a_{34}a_{45}\overline{a_{15}}.$

\item $r^2[\alpha_w\alpha_t\alpha_v - (\alpha_w\mu_t\mu_v   +\alpha_t\mu_w\mu_v  
+\alpha_v\mu_w\mu_t ) + 2\mu_w\mu_t\mu_v \cos(\theta)]\\
+4\alpha_p\alpha_q [(\alpha_w\mu_t\mu_v   +\alpha_t\mu_w\mu_v  +\alpha_v\mu_w\mu_t ) 
- 2\mu_w\mu_t\mu_v \cos(\theta)]\\
=\sum\limits_{S_{ijk}S_{lm}}|a_{lm}|^2 \alpha_i\alpha_j\alpha_k 
-\sum\limits_{S_{ij}S_{klm}}\Re(a_{kl}a_{lm}\overline{a_{km}})\alpha_i\alpha_j 
-\dfrac{1}{4}\sum\limits_{i=1}^5P_i\alpha_i 
+\dfrac{1}{2}\sum\limits_{S_{i}S_{jklm}}\Re(a_{jk}a_{kl}a_{lm}\overline{a_{jm}})\alpha_i\\
+\dfrac{1}{4}\sum\limits_{S_{ij}S_{klm}}\Re(a_{kl}a_{lm}\overline{a_{km}}) |a_{ij}|^2 
-\dfrac{1}{4}\sum\limits_{S_{ijkl}S_{iml}}\Re(a_{ij}a_{jk}a_{kl}\overline{a_{im}}\overline{a_{ml}})\\
-\dfrac{1}{4}\sum\limits_{S_{ijk}S_{lm}S_{im}S_{lk}}\Re(a_{ij}a_{jk}a_{lm}\overline{a_{im}}\overline{a_{lk}})
-\dfrac{1}{4}\Re(a_{12}a_{23}a_{34}a_{45}\overline{a_{15}}).$

\item $r^2[\beta_w\beta_t\beta_v - (\beta_w\mu_t\mu_v   +\beta_t\mu_w\mu_v  +\beta_v\mu_w\mu_t )
 + 2\mu_w\mu_t\mu_v \sin(\theta)]\\
+4\beta_p\beta_q [(\beta_w\mu_t\mu_v   +\beta_t\mu_w\mu_v  +\beta_v\mu_w\mu_t )
 - 2\mu_w\mu_t\mu_v \sin(\theta)]\\
 =\sum\limits_{S_{ijk}S_{lm}}|a_{lm}|^2 \beta_i\beta_j\beta_k 
-\sum\limits_{S_{ij}S_{klm}}\Im(a_{kl}a_{lm}\overline{a_{km}})\beta_i\beta_j 
-\dfrac{1}{4}\sum\limits_{i=1}^5P_i\beta_i 
-\dfrac{1}{2}\sum\limits_{S_{i}S_{jklm}}\Re(a_{jk}a_{kl}a_{lm}\overline{a_{jm}})\beta_i\\
+\dfrac{1}{4}\sum\limits_{S_{ij}S_{klm}}\Im(a_{kl}a_{lm}\overline{a_{km}}) |a_{ij}|^2 
-\dfrac{1}{4}\sum\limits_{S_{ijkl}S_{iml}}\Im(a_{ij}a_{jk}a_{kl}\overline{a_{im}}\overline{a_{ml}})\\
-\dfrac{1}{4}\sum\limits_{S_{ijk}S_{lm}S_{im}S_{lk}}\Im(a_{ij}a_{jk}a_{lm}\overline{a_{im}}\overline{a_{lk}})
+\dfrac{1}{4}\Im(a_{12}a_{23}a_{34}a_{45}\overline{a_{15}}).$

 \item  $r^2 [ (\alpha_t\alpha_v   + \alpha_w\alpha_v  +\alpha_w\alpha_t )-(\mu_t\mu_v   + \mu_w\mu_v  +\mu_w\mu_t )] +4\alpha_p\alpha_q(\mu_t\mu_v   + \mu_w\mu_v  +\mu_w\mu_t )\\
 +4(\alpha_p+\alpha_q)[(\alpha_w\mu_t\mu_v   +\alpha_t\mu_w\mu_v  +\alpha_v\mu_w\mu_t ) - 2\mu_w\mu_t\mu_v \cos(\theta)]\\
 =  \sum\limits_{S_{ijk}S_{lm}}|a_{lm}|^2 (\alpha_i\alpha_j+\alpha_i\alpha_k + \alpha_j\alpha_k)
- \sum\limits_{S_{ij}S_{klm}}\Re(a_{kl}a_{lm}\overline{a_{km}})(\alpha_i+\alpha_j)\\
-\dfrac{1}{4}\sum\limits_{i=1}^5P_i
+ \dfrac{1}{2}\sum\limits_{S_{jklm}}\Re(a_{jk}a_{kl}a_{lm}\overline{a_{jm}}).
 $ 
\item $\mu_w\mu_t\mu_v \neq 0$
\item  $\lambda_j\mu_k + \lambda_k\mu_j - 2\mu_j\mu_k e^{i\theta}\neq \lambda_i(\mu_k + \mu_j)$
for every $i,j,k \in \lbrace w,t,v \rbrace: i\neq j\neq k. $
\end{enumerate}
 where 
$P_i=\sum\limits_{S_{jkm}S_{lm}}|a_{jk}|^2|a_{lm}|^2-\sum\limits_{S_{jkl}S_{jml}}\Re(a_{jk}a_{kl}\overline{a_{jm}}\overline{a_{ml}})-\sum\limits_{S_{jk}S_{lkm}}\Re(a_{jk}a_{lm}\overline{a_{jm}}\overline{a_{lk}}).$\\ for every $i=1,..,5$. 
\end{theorem}
%
%
\begin{proof}
Taking conditions from lemma \ref{flat}
 \item    \begin{align*} 
 P_A(x, y, z)& =\left[(\alpha_px+\beta_py+z)(\alpha_qx+\beta_qy+z)-\frac{r^2}{4}(x^2 +y^2)
  \right]\\ & [(\alpha_wx+\beta_wy+z)(\alpha_tx+\beta_ty+z)(\alpha_vx+\beta_vy+z)\\
   & -(\alpha_wx+\beta_wy+z)(x^2 +y^2)\mu_t\mu_v\\ & 
  -(\alpha_tx+\beta_ty+z)(x^2 +y^2)\mu_w\mu_v  \\ & 
  -(\alpha_vx+\beta_vy+z)(x^2 +y^2)\mu_w\mu_t  \\
   &+(x-iy)(x^2 +y^2)\mu_w\mu_t\mu_ve^{i\theta} \\
  & + (x+iy)(x^2 +y^2)\mu_w\mu_t\mu_ve^{-i\theta}
    \end{align*}
\item $\mu_w\mu_t\mu_v \neq 0$
\item  $\lambda_j\mu_k + \lambda_k\mu_j -2 \mu_j\mu_k e^{i\theta}\neq \lambda_i(\mu_k + \mu_j)$
for every $i,j,k \in \lbrace w,t,v \rbrace: i\neq j\neq k.$\\
Comparing the previous formula of $P_A(x, y, z)$ with polynomial $(\ast)$. We find 
  \begin{align*} 
 Q(x, y, z)& =r^2[(\alpha_wx+\beta_wy+z)(\alpha_tx+\beta_ty+z)(\alpha_vx+\beta_vy+z)\\
   & - (x^2 +y^2)\left((\alpha_wx+\beta_wy+z)\mu_t\mu_v  + (\alpha_tx+\beta_ty+z)\mu_w\mu_v   
 +(\alpha_vx+\beta_vy+z)\mu_w\mu_t\right) \\
 &+2(x^2 +y^2)(x\mu_w\mu_t\mu_v \cos(\theta)+y\mu_w\mu_t\mu_v \sin(\theta))] \\
 & + 4(\alpha_px+\beta_py+z)(\alpha_qx+\beta_qy+z)\\
 & [\left((\alpha_wx+\beta_wy+z)\mu_t\mu_v  + (\alpha_tx+\beta_ty+z)\mu_w\mu_v   
 +(\alpha_vx+\beta_vy+z)\mu_w\mu_t\right) \\
 &-2(x\mu_w\mu_t\mu_v \cos(\theta)+y\mu_w\mu_t\mu_v \sin(\theta))].
    \end{align*}
    
 Computing the coefficients of $x^3,y^3,z^3,x^2y,xy^2,x^2z,xz^2,y^2z,yz^2,xyz$ by identification, we find respectively 
\begin{enumerate}

\item $r^2[\alpha_w\alpha_t\alpha_v - (\alpha_w\mu_t\mu_v   +\alpha_t\mu_w\mu_v  
+\alpha_v\mu_w\mu_t ) + 2\mu_w\mu_t\mu_v \cos(\theta)]\\
+4\alpha_p\alpha_q [(\alpha_w\mu_t\mu_v   +\alpha_t\mu_w\mu_v  +\alpha_v\mu_w\mu_t ) 
- 2\mu_w\mu_t\mu_v \cos(\theta)]\\
=\sum\limits_{S_{ijk}S_{lm}}|a_{lm}|^2 \alpha_i\alpha_j\alpha_k 
-\sum\limits_{S_{ij}S_{klm}}\Re(a_{kl}a_{lm}\overline{a_{km}})\alpha_i\alpha_j 
-\dfrac{1}{4}\sum\limits_{i=1}^5P_i\alpha_i 
+\dfrac{1}{2}\sum\limits_{S_{i}S_{jklm}}\Re(a_{jk}a_{kl}a_{lm}\overline{a_{jm}})\alpha_i\\
+\dfrac{1}{4}\sum\limits_{S_{ij}S_{klm}}\Re(a_{kl}a_{lm}\overline{a_{km}}) |a_{ij}|^2 
-\dfrac{1}{4}\sum\limits_{S_{ijkl}S_{iml}}\Re(a_{ij}a_{jk}a_{kl}\overline{a_{im}}\overline{a_{ml}})\\
-\dfrac{1}{4}\sum\limits_{S_{ijk}S_{lm}S_{im}S_{lk}}\Re(a_{ij}a_{jk}a_{lm}\overline{a_{im}}\overline{a_{lk}})
-\dfrac{1}{4}\Re(a_{12}a_{23}a_{34}a_{45}\overline{a_{15}}).$

\item $r^2[\beta_w\beta_t\beta_v - (\beta_w\mu_t\mu_v   +\beta_t\mu_w\mu_v  +\beta_v\mu_w\mu_t )
 + 2\mu_w\mu_t\mu_v \sin(\theta)]\\
+4\beta_p\beta_q [(\beta_w\mu_t\mu_v   +\beta_t\mu_w\mu_v  +\beta_v\mu_w\mu_t )
 - 2\mu_w\mu_t\mu_v \sin(\theta)]\\
 =\sum\limits_{S_{ijk}S_{lm}}|a_{lm}|^2 \beta_i\beta_j\beta_k 
-\sum\limits_{S_{ij}S_{klm}}\Im(a_{kl}a_{lm}\overline{a_{km}})\beta_i\beta_j 
-\dfrac{1}{4}\sum\limits_{i=1}^5P_i\beta_i 
-\dfrac{1}{2}\sum\limits_{S_{i}S_{jklm}}\Re(a_{jk}a_{kl}a_{lm}\overline{a_{jm}})\beta_i\\
+\dfrac{1}{4}\sum\limits_{S_{ij}S_{klm}}\Im(a_{kl}a_{lm}\overline{a_{km}}) |a_{ij}|^2 
-\dfrac{1}{4}\sum\limits_{S_{ijkl}S_{iml}}\Im(a_{ij}a_{jk}a_{kl}\overline{a_{im}}\overline{a_{ml}})\\
-\dfrac{1}{4}\sum\limits_{S_{ijk}S_{lm}S_{im}S_{lk}}\Im(a_{ij}a_{jk}a_{lm}\overline{a_{im}}\overline{a_{lk}})
+\dfrac{1}{4}\Im(a_{12}a_{23}a_{34}a_{45}\overline{a_{15}}).$

 \item  $r^2 +4(\mu_t\mu_v   + \mu_w\mu_v  +\mu_w\mu_t )= \sum\limits_{S_{lm}}|a_{lm}|^2.$

 \item  $r^2 [(\beta_w\alpha_t\alpha_v   +\beta_t\alpha_w\alpha_v  +\beta_v\alpha_w\alpha_t ) 
 + 2\mu_w\mu_t\mu_v \sin(\theta)- (\beta_w\mu_t\mu_v   +\beta_t\mu_w\mu_v  +\beta_v\mu_w\mu_t ) ] \\
 +4\alpha_p\alpha_q [(\beta_w\mu_t\mu_v   +\beta_t\mu_w\mu_v  +\beta_v\mu_w\mu_t ) - 2\mu_w\mu_t\mu_v \sin(\theta)] \\
 +4(\alpha_p\beta_q + \alpha_q\beta_p)[(\alpha_w\mu_t\mu_v   +\alpha_t\mu_w\mu_v  +\alpha_v\mu_w\mu_t )- 2\mu_w\mu_t\mu_v \cos(\theta)]\\
 =\sum\limits_{S_{ijk}S_{lm}}|a_{lm}|^2 (\alpha_i\alpha_j\beta_k+\alpha_i\alpha_k\beta_j + \alpha_j\alpha_k\beta_i)
 -\sum\limits_{S_{ij}S_{klm}}\Im(a_{kl}a_{lm}\overline{a_{km}})\alpha_i\alpha_j \\
 -\sum\limits_{S_{ij}S_{klm}}\Re(a_{kl}a_{lm}\overline{a_{km}})(\alpha_i\beta_j +\alpha_j\beta_i) 
 -\dfrac{1}{4}\sum\limits_{i=1}^5P_i\beta_i 
 +\dfrac{1}{2}\sum\limits_{S_{i}S_{jklm}}\Re(a_{jk}a_{kl}a_{lm}\overline{a_{jm}})\beta_i\\
 +\sum\limits_{S_{i}S_{jklm}}\Im(a_{jk}a_{kl}a_{lm}\overline{a_{jm}})\alpha_i
 + \dfrac{1}{4}\sum\limits_{S_{ij}S_{klm}}\Im(a_{kl}a_{lm}\overline{a_{km}}) |a_{ij}|^2\\
 -\dfrac{1}{4}\sum\limits_{S_{ijkl}S_{iml}}\Im(a_{ij}a_{jk}a_{kl}\overline{a_{im}}\overline{a_{ml}})
 -\dfrac{1}{4}\sum\limits_{S_{ijk}S_{lm}S_{im}S_{lk}}\Im(a_{ij}a_{jk}a_{lm}\overline{a_{im}}\overline{a_{lk}}) \\
-\dfrac{3}{4}\Im(a_{12}a_{23}a_{34}a_{45}\overline{a_{15}})  
. $

  \item  $r^2 [(\alpha_w\beta_t\beta_v   +\alpha_t\beta_w\beta_v  +\alpha_v\beta_w\beta_t ) 
 + 2\mu_w\mu_t\mu_v \cos(\theta)- (\alpha_w\mu_t\mu_v   +\alpha_t\mu_w\mu_v  +\alpha_v\mu_w\mu_t ) ] \\
 +4\beta_p\beta_q [(\alpha_w\mu_t\mu_v   +\alpha_t\mu_w\mu_v  +\alpha_v\mu_w\mu_t ) - 2\mu_w\mu_t\mu_v \cos(\theta)] \\
 +4(\alpha_p\beta_q + \alpha_q\beta_p)[(\beta_w\mu_t\mu_v   +\beta_t\mu_w\mu_v  +\beta_v\mu_w\mu_t ) - 2\mu_w\mu_t\mu_v \sin(\theta)]\\
 =\sum\limits_{S_{ijk}S_{lm}}|a_{lm}|^2 (\alpha_i\beta_j\beta_k+\alpha_j\beta_i\beta_k + \alpha_k\beta_i\beta_j)
 -\sum\limits_{S_{ij}S_{klm}}\Re(a_{kl}a_{lm}\overline{a_{km}})\beta_i\beta_j \\
 -\sum\limits_{S_{ij}S_{klm}}\Im(a_{kl}a_{lm}\overline{a_{km}})(\alpha_i\beta_j +\alpha_j\beta_i)  
 -\dfrac{1}{4}\sum\limits_{i=1}^5P_i\alpha_i  
 -\dfrac{1}{2}\sum\limits_{S_{i}S_{jklm}}\Re(a_{jk}a_{kl}a_{lm}\overline{a_{jm}})\alpha_i\\ 
 +\sum\limits_{S_{i}S_{jklm}}\Im(a_{jk}a_{kl}a_{lm}\overline{a_{jm}})\beta_i 
 + \dfrac{1}{4}\sum\limits_{S_{ij}S_{klm}}\Re(a_{kl}a_{lm}\overline{a_{km}}) |a_{ij}|^2\\ 
 -\dfrac{1}{4}\sum\limits_{S_{ijkl}S_{iml}}\Re(a_{ij}a_{jk}a_{kl}\overline{a_{im}}\overline{a_{ml}}) 
 -\dfrac{1}{4}\sum\limits_{S_{ijk}S_{lm}S_{im}S_{lk}}\Re(a_{ij}a_{jk}a_{lm}\overline{a_{im}}\overline{a_{lk}}) \\
+\dfrac{3}{4}\Re(a_{12}a_{23}a_{34}a_{45}\overline{a_{15}}) .  $

 \item  $r^2 [ (\alpha_t\alpha_v   + \alpha_w\alpha_v  +\alpha_w\alpha_t )-(\mu_t\mu_v   + \mu_w\mu_v  +\mu_w\mu_t )] +4\alpha_p\alpha_q(\mu_t\mu_v   + \mu_w\mu_v  +\mu_w\mu_t )\\
 +4(\alpha_p+\alpha_q)[(\alpha_w\mu_t\mu_v   +\alpha_t\mu_w\mu_v  +\alpha_v\mu_w\mu_t ) - 2\mu_w\mu_t\mu_v \cos(\theta)]\\
 =  \sum\limits_{S_{ijk}S_{lm}}|a_{lm}|^2 (\alpha_i\alpha_j+\alpha_i\alpha_k + \alpha_j\alpha_k)
- \sum\limits_{S_{ij}S_{klm}}\Re(a_{kl}a_{lm}\overline{a_{km}})(\alpha_i+\alpha_j)\\
-\dfrac{1}{4}\sum\limits_{i=1}^5P_i
+ \dfrac{1}{2}\sum\limits_{S_{jklm}}\Re(a_{jk}a_{kl}a_{lm}\overline{a_{jm}}).
 $

 \item $r^2(\alpha_w+\alpha_t+\alpha_v)+4(\alpha_w\mu_t\mu_v   +\alpha_t\mu_w\mu_v  +\alpha_v\mu_w\mu_t ) - 8\mu_w\mu_t\mu_v \cos(\theta)\\ +4(\alpha_p+\alpha_q)(\mu_t\mu_v   + \mu_w\mu_v  +\mu_w\mu_t )\\
 =\sum\limits_{S_{ijk}S_{lm}}|a_{lm}|^2 (\alpha_i + \alpha_j + \alpha_k) 
 -\sum\limits_{S_{klm}}\Re(a_{kl}a_{lm}\overline{a_{km}}).
 $

  \item  $r^2 [ (\beta_t\beta_v   + \beta_w\beta_v  +\beta_w\beta_t )-(\mu_t\mu_v   +  \mu_w\mu_v  +\mu_w\mu_t )] +4\beta_p\beta_q(\mu_t\mu_v   + \mu_w\mu_v  +\mu_w\mu_t )\\
 +4(\beta_p+\beta_q)[(\beta_w\mu_t\mu_v   +\beta_t\mu_w\mu_v  +\beta_v\mu_w\mu_t ) - 2\mu_w\mu_t\mu_v \sin(\theta)] \\
 =    \sum\limits_{S_{ijk}S_{lm}}|a_{lm}|^2 (\beta_i\beta_j+\beta_i\beta_k + \beta_j\beta_k)
- \sum\limits_{S_{ij}S_{klm}}\Im(a_{kl}a_{lm}\overline{a_{km}})(\beta_i+\beta_j)\\
-\dfrac{1}{4}\sum\limits_{i=1}^5P_i
-\dfrac{1}{2}\sum\limits_{S_{jklm}}\Re(a_{jk}a_{kl}a_{lm}\overline{a_{jm}}).$  
 \item $r^2(\beta_w+\beta_t+\beta_v)+4(\beta_w\mu_t\mu_v   +\beta_t\mu_w\mu_v  +\beta_v\mu_w\mu_t ) - 8\mu_w\mu_t\mu_v \sin(\theta)\\ +4(\beta_p+\beta_q)(\mu_t\mu_v   + \mu_w\mu_v  +\mu_w\mu_t )\\
 =  \sum\limits_{S_{ijk}S_{lm}}|a_{lm}|^2 (\beta_i + \beta_j + \beta_k) 
 -\sum\limits_{S_{klm}}\Im(a_{kl}a_{lm}\overline{a_{km}}).
 $
 \item 
 $r^2(\alpha_w\beta_t +\alpha_w\beta_v+\alpha_t\beta_w+\alpha_t\beta_v+\alpha_v\beta_w+\alpha_v\beta_t)+ 4(\alpha_p\beta_q+\alpha_q\beta_p)(\mu_t\mu_v   + \mu_w\mu_v  +\mu_w\mu_t )\\
 +4(\alpha_p+\alpha_q)[(\beta_w\mu_t\mu_v   +\beta_t\mu_w\mu_v  +\beta_v\mu_w\mu_t ) - 2\mu_w\mu_t\mu_v \sin(\theta)] \\
  +4(\beta_p+\beta_q)[(\alpha_w\mu_t\mu_v   +\alpha_t\mu_w\mu_v  +\alpha_v\mu_w\mu_t ) - 2\mu_w\mu_t\mu_v \cos(\theta)]\\
 = \sum\limits_{S_{ijk}S_{lm}}|a_{lm}|^2 (\alpha_i\beta_j +\alpha_j\beta_i+\alpha_i\beta_k+\alpha_k\beta_i+\alpha_j\beta_k+\alpha_k\beta_j)- \sum\limits_{S_{ij}S_{klm}}\Re(a_{kl}a_{lm}\overline{a_{km}})(\beta_i+\beta_j) \\
  - \sum\limits_{S_{ij}S_{klm}}\Im(a_{kl}a_{lm}\overline{a_{km}})(\alpha_i+\alpha_j) 
+\sum\limits_{S_{jklm}}\Im(a_{jk}a_{kl}a_{lm}\overline{a_{jm}})$.\\

\end{enumerate}
Note that the combination of $(1)$, $(2)$, $(4)$ and $(5)$ is equivalent to( $d)$, $(e)$ and $(f)$, since 
$(1)-(5)-i(2)+i(4)$ yields (d).

 $(6)$, $(8)$ and $(10)$ is
equivalent to $(c)$ and $(g)$ as $(6)-(8) + i(10)$ yields $(c)$.

$(7)$ and $(9)$ is equivalent to $(b)$ as it's follows from $(7)+i(9)$.

 This completes the proof.
\end{proof} 
\section{On the circular numerical range of $S_5$ matrices }

       Recall that an n-by-n
       matrix $A$ is said to be of class $S_n$ if $A$ is a contraction, the eigenvalues 
       of $A$ are all in the open unit disc $\mathbb{D}$ and
       $rank (I_n -A^*A) = 1$.
Two unitary equivalent $S_n$-matrices have the following useful characterization. 
\begin{lemma}\cite[Theorem 4.1]{GauSn}\label{lm2}
$A_1$ and $A_2$ be n-dimensional operators with $A_2$ in $S_n$.
Then $A_1$ is unitary equivalent to $A_2$ if and only if $A_1$ is a contraction and $W(A_1)=W(A_2)$.
\end{lemma}

Now, we establish the result of Theorem \ref{5by5iso} in the special case of $S_5$-matrices. 
\begin{theorem}\label{S_5}
Let $A$ be  a $S_5$ matrix with $W(A)=\left\lbrace z\in \mathbb{C}:  \left|z-a \right|\leq r \right\rbrace$, ($r>0)$.
if kippenhahn curve $C_R(A)$ have one of the following   shape
\begin{enumerate}[label=(\roman*)]
\item $C_R(A)$ consist of three point and an elliptic disc.
\item $C_R(A)$ consist of two elliptic disc and a point. 
\end{enumerate}
 then $a=0$.
\end{theorem}

\begin{proof}
 Without
loss of generality, we assume that  $A$ in an upper triangular matrix. The assumption on the numerical range of $ A$ imply  that the origin  $a$ is an eigenvalue  with algebraic multiplicity at least 2.  So, by \cite[Corollary 1.3]{GauLuc} $A$  takes the following form
$$A=\left[ \begin {array}{ccccc} a&1-{a}^{2}&-a\sqrt {1-{a}^{2}}&0&0
\\ \noalign{\medskip}0&a&\sqrt {1-{a}^{2}}&0&0\\ \noalign{\medskip}0&0
&0&\sqrt {1-  \left| b \right|^{2}}&-\overline{b}\sqrt {1-
   \left| c \right|  ^{2}}\\ \noalign{\medskip}0&0&0&b&
\sqrt {1- \left| b \right|  ^{2}}\sqrt {1-
 \left| c \right|   ^{2}}\\ \noalign{\medskip}0&0&0&0&c
\end {array} \right]. $$
Moreover, we can take $a$ positive by a suitable rotation, thus $W(A)$  is symmetric
       with respect to the real axis, witch means  that $W(A) = W(A^*)$, ($A^*$ is the adjoint matrix of $A$),
       as we mentioned  below $A$ is $S_5$ 
       and therefore by Lemma \ref{lm2} $A$ and $A^*$ are unitary equivalent, 
       moreover one can see that the eigenvalues $b$ and $c$ of $A$ must  be real or complex conjugates . Let
$$B= A-aI_5=\left[ \begin {array}{ccccc} 0&1-{a}^{2}&-a\sqrt {1-{a}^{2}}&0&0
\\ \noalign{\medskip}0&0&\sqrt {1-{a}^{2}}&0&0\\ \noalign{\medskip}0&0
&-a&\sqrt {1-  |b|^{2}}&-\overline{b}\sqrt {1-  |c|^{2}}\\ \noalign{\medskip}0&0&0&b-a&
\sqrt {1-   |b|^{2}}\sqrt {1-|c|^{2}}\\ \noalign{\medskip}0&0&0&0&c-a
\end {array} \right].$$ 
        Consider the homogeneous  Kippenhahn polynomial $p_B(x, y, z) =det(x\Re B + y\Im B + zI_5)$ 
        of degree $5$ on the complex projective plane $\mathbb{CP}^2$. Since by
         hypotheses  $W(B)$ is a circular disc with center $0$ and radius $r$, By 
         assumption $C_R(B)$ have one of two possibles shapes.
\begin{itemize}
	\item [(i)] A circle with center $0$ and radius $r$ together 
        with  three points $-a,b-a,c-a$ inside it.	
        \item [(ii)] A circle with center $0$ and radius $r$ together 
        with a point($-a,b-a$ or $c-a$) and an ellipse with (minor axis length $s\leq 2r $) and the two remaining points as the foci, all inside the circle.
\end{itemize}
 Applying condition $(d)$ of Theorem \ref{th2ellipse} to the upper-triangular matrix $B$
         yields to 
\begin{align*}
4r^2 &(-a(b-a)(c-a)) \\
&= (1-a^2 )^2 (-a)(b-a)(c-a)-(1-a^2 )(-a) \sqrt{1-a^2 } \sqrt {1- a^2 }(b-a)(c-a)\\
&=0.
\end{align*}
                  Then either $a = 0$, $a = b$ or $a = c$.  If it is first so 
                  it will done otherwise if it is one of the two  latter, the
                   condition $(c)$  of Theorem \ref{th2ellipse}  gives,
\begin{align*}            
4r^2&(-a(b-a)+(b-a)(c-a)-a(c-a))\\
&=(1-a^2 )^2(-a(b-a)+(b-a)(c-a)-a(c-a))\\
&+a^2(1-a^2)(b-a)(c-a)+(1-a^2)(b-a)(c-a)\\
&-(1-a^2)^2(-a(b-a)-a(c-a))\\
&=0.
\end{align*}
                   Thus $a=b$ if $a=c$ and vise versa. By the condition $(b)$ of 
                   Theorem \ref{th2ellipse}
 \begin{align*}                 
4r^2&(-a+(b-a)+(c-a))+s^2(0+0+\lambda)= (1-a^2)^2(-a+b-a+c-a)
\\&+a^2(1-a^2)(b-a+c-a) + (1-a^2)(b-a+c-a)+(1-b^2)(c-a)+b^2(1-c^2)(b-a) \\
&-a(1-b^2)(1-c^2) +a(1-a^2)^2+ b(1-b^2)(1-c^2).
\end{align*}
 where $\lambda$ take one of  the eigenvalue $-a, b-a$ or $c-a$.\\           
                   Assuming that $a=b=c$, we
                   fined  $(4r^2+s^2)(-a)=0$ or $ar^2=0$. Hence $a=0$.\\ 
This complete the proof.
\end{proof}
%
%
\begin{remark}\label{RkSn}
As it well known that for every  $S_n$ matrix $T$, $\Re(T)$ have only simple eigenvalues see \cite[Corollary 2.6]{GauSn}, then $C_R(A)$  may not contain a curve of degree 4 with double tangent.
\end{remark}
\section{Proof of Theorem \ref{5by5iso}}
 
%
%
%
It well known that a n-by-n partial isometry $A$ can be represented $Ker(A)\oplus Ker(A)^{\perp}$, by
$$A=\begin{bmatrix}0&B
\\ 0&C
\end{bmatrix} $$
with $B$ and $C$ satisfying $B^*B+C^*C=I$, see \cite[Proposition 2.1]{GauPI}. 
Another result which we need for the proof of Theorem \ref{5by5iso}  concerns the irreducibility
of a partial isometry.
%
%
%
\begin{lemma}\cite[lemma 2.8]{GauPI}\label{lm3} 
    Let 
$$A=\begin{bmatrix}
0_m&B
\\ 0& C
\end{bmatrix} \quad \text{ on }\quad \mathbb{C}^n =\mathbb{C}^m\oplus \mathbb{C}^{n-m}, \quad (1 \leq m \leq n).
$$
\begin{itemize}
	\item[(a)] If $k = rank B < m$, then $A$ is unitarily similar to $0_{m-k} \oplus$ $A_1$ for some matrix
	$$A_1=\begin{bmatrix}
	0_k&B_1
	\\ 0& C_1
	\end{bmatrix} \quad \text{ on }\quad \mathbb{C}^{n-m+k} =\mathbb{C}^k\oplus \mathbb{C}^{n-m}, \quad \text{ with }\quad rank B_1 = k.$$
	\item[(b)] If $m > [n/2]$, the largest integer less than or
	equal to $n/2$, then A is  unitarily similar to $0_{2m-n} \oplus A_2$ for some matrix 
	$$A_2=\begin{bmatrix}
	0_{n-m}&B_2
	\\ 0& C_2
	\end{bmatrix} \quad \text{ on }\quad \mathbb{C}^{2(n-m)} =\mathbb{C}^{n-m}\oplus \mathbb{C}^{n-m}.$$
\end{itemize}
\end{lemma}
The next proposition relates partial
isometries with $S_n$-matrices.
\begin{proposition}\cite[Proposition 2.3]{GauPI}\label{lm1} 
	Let $A$ be an n-by-n matrix. Then $A$ is an irreducible partial isometry
	with $dim ker A = 1$ if and only if $A$ is of class $S_n$ with 0 in $\sigma(A)$.
\end{proposition}

Now, we are ready to establish our main theorem.
\begin{proof}[Proof of Theorem \ref{5by5iso}]
Let $A$ be an $5\times5$ partial isometry with 
           $W(A)=\left\{ z\in \mathbb{C}:  \left|z-a \right|\leq r \right\}$, ($r>0$). First 
           let us remark that if $A$ is reducible, then $A$ is unitarily 
           similar to $A_1\oplus A_2$, where  $A_1$ and $A_2$ are two partial 
           isometry with order at most 4.\\ 
 Since  one of $W(A_1)$ or $W(A_2)$ must be equal to that of $A$ , so by Theorem \ref{thm4by4}
          it follows that $a=0$.
          
  Now, we assume that $A$ is irreducible. According to the dimension of the kernel of $A$, we
        distinguish three cases
        
\textbf{Case 1.} $\dim ker A = 1$. By the Proposition \ref{lm1}, $A$ is $S_5$-matrix, so according to Theorem \ref{S_5} and Remark \ref{RkSn}, $a=0$.\\
%
%
\textbf{Case 2.}  
$dim ker A = 2$. As $W(A)$ is a circular disc centered at $a$, we may assume that
$$A= \begin{bmatrix} 0&B \\ \noalign{\medskip} 0&C\end{bmatrix}= \begin{bmatrix} 0&0&k&l&t\\ 
0&0&g&h&j
\\ 0&0&b&e&f\\ 
0&0&0&a&d
\\0&0&0&0&a\end{bmatrix} \quad \text{ on } \mathbb{C}^2\oplus\mathbb{C}^3,$$  
 with
 \begin{align*}
I_3&= B^*B+C^*C\\
&=\begin{bmatrix} \left| k \right|^{2}+
 \left| g \right|^{2}+  \left| b \right|^{2}&\overline{k}l+\overline{g}h+\overline{b}e&\overline{k}t+
\overline{g}j+\overline{b}f\\ \noalign{\medskip}
\overline{l}k+
\overline{h}g+\overline{e}b&  \left| l \right|^{2}+
\left| h \right|^{2}+  \left| e \right|^{2}+ \left| a \right|^{2}&\overline{l}t+
\overline{h}j+\overline{e}f+\overline{a}d\\ 
\overline{t}k+\overline{j}g+\overline{f}b&\overline{t}l+\overline{j}h+
\overline{f}e+\overline{d}a&  \left| t \right|^{2}+
\left| j \right|^{2}+  \left| f \right|^{2}+ \left| d \right|^{2}+ \left|a \right|^{2}\end{bmatrix}.
 \end{align*}
As in the proof of Theorem\ref{S_5}, $a$ is positive and  $C_R(A)$ have one of the three possible shapes.

\begin{itemize}
	\item [(i)]$C_R(A)$ contains a disc (with centre $a$, radius $r$) and  three point $0,0,b$.
	\item [(ii)]$C_R(A)$ is a disc (with centre $a$, radius $r$) , together with  and an ellipse and a point.
	\item [(iii)] $C_R(A)$ contains a disc (with centre $a$, radius $r$) and a 4 degree  curve with a double tangent.
\end{itemize}
  Applying condition $(d)$ of Theorem \ref{th2ellipse} and Theorem \ref{th4degree} to $A-aI_5$ we get 
\begin{align*}
4r^2(a.a.(b-a))&=|d|^2a^2(b-a)-a^2ed\overline{f}+a(b-a)hd\overline{j}+a(b-a)ld\overline{t}\\
&-aged\overline{j}-aked\overline{t}\\
&=|d|^2a^2(b-a)-a^2ed\overline{f}+ad(b-a)(h\overline{j}+l\overline{t}) -aed(g\overline{j} +k\overline{t})\\
&=|d|^2a^2(b-a)-a^2ed\overline{f}-ad(b-a)(e\overline{f}+a\overline{d})+aedb\overline{f}\\
&=0
\end{align*}
thus $a=0$ or $b=a$ if it'is the later, by condition $(c)$ of Theorem \ref{th2ellipse} and Theorem  \ref{th4degree}
\begin{align*}
4r^2(a^2)& = a^2(|e|^2+|f|^2+|d|^2) +2aed\overline{f} +aeg\overline{h}+aek\overline{l} + afg\overline{j} + afk\overline{t} \\
&+ adh\overline{j}+adl\overline{t}+edg\overline{j} + edk\overline{t}.\\
&= a^2(|e|^2+|f|^2+|d|^2) +2aed\overline{f} + ae(g\overline{h} + k\overline{l})+af(g\overline{j} + k\overline{t})\\
&+ad(h\overline{j}+l\overline{t})+ed(g\overline{j}+k\overline{t})\\
&=0
\end{align*}
and therefore $a=0$.

\textbf{Case 3.}  
$dim ker A > 2$, then it follow from  Lemma \ref{lm3} that $A$ is reducible, then $a=0$. 

 This completes the proof of the theorem.
\end{proof}

\begin{remark} To give a  complete answer to the conjecture of Gau et al,     in dimension 5, it remains to study  the case when $C_R(A)$ is  an elliptic disc and a curve of order 6, consisting of an oval and a curve of three cups.
	
Based on the factoribility of $p_A$, Kippenhahn  in \cite{Kippenhahn1951} gives a fully classification   of  the numerical range of $3\times 3$ matrices, and a pertinent tests were offered in \cite{Rodman}, however there is no much provided about the connection between concrete description of the curve $C_R(A)$ and $p_A$  when $W(A)$ is oval, as like  given in Lemma \ref{flat} when the boundary of $W(A)$ have a flat portion. Thus, the remaining case is a difficult question to treated.
\end{remark}


\end{document}